\documentclass[aps,pra,preprint,superscriptaddress]{revtex4}
\usepackage{amsmath,amssymb,amsfonts}
\usepackage[english]{babel}
\usepackage{graphicx}
\usepackage{dcolumn}
\usepackage{bm}
\usepackage{color}

\begin{document} 

\title{Quantifying momenta through the Fourier transform} 

\author{B.~M. Rodr\'{\i}guez-Lara}
\email[]{cqtbmrl@nus.edu.sg} 
\affiliation{Centre for Quantum Technologies, National University of Singapore, 2 Science Drive 3, Singapore 117542.}

\begin{abstract} Integral transforms arising from the separable solutions to the Helmholtz differential equation are presented. Pairs of these integral transforms are related via Plancherel theorem and, ultimately, any of these integral transforms may be calculated using only Fourier transforms. This result is used to evaluate the mean value of momenta associated to the symmetries of the reduced wave equation.
As an explicit example, the orbital angular momenta of plane and elliptic-cylindrical waves is presented. 
\end{abstract} 

\maketitle

\section{Introduction} 
The orbital angular momentum (OAM) of light has become a fundamental concept in optics \cite{Allen2003,Andrews2008}. The fact that light may carry orbital angular momentum which is not related to its polarization state has been shown both theoretically \cite{Beth1936p115} and experimentally \cite{Allen1992p8185}. The optical tweezers toolbox now includes experiments involving the transfer of OAM from paraxial beams \cite{He1995p826, Milne2007p13972} and vector fields \cite{Wang2010p253602} to dielectric microscopic particles; also, the feasibility of transferring OAM from laser beams into cold atoms \cite{Tabosa1999p4967,Picon10p083053}, and even store it in matter ensembles \cite{Moretti2009p023825}, has been discussed and shown.
Furthermore, a couple of dynamical variables of light, involving compositions of linear and angular momenta, have been theoretically identified related to propagation invariant vector fields with transversal elliptic \cite{RodriguezLara2008p033813} and parabolic \cite{RodriguezLara2009p055806} symmetries. In the paraxial regime, the elliptic \cite{GutierrezVega2000p1493} and parabolic \cite{Bandres2004p44} beams are well known; they have been generated by diverse means \cite{ChavezCerda2002S52,AlvarezElizondo2008p18770,Hernandez2010p366903} and the orbital angular momentum per photon of the elliptic beams discussed \cite{LopezMariscal2005p2364}. 
Proposed applications for these light fields carrying compositions of linear and angular momenta include the optical manipulation of atoms, \textit{e.g.} thermal atomic cloud splitters \cite{RodriguezLara2009p011813R} and atomic billiards \cite{PerezPascual2011p}, as well as the optical manipulation of dielectric microscopic particles in elliptic, hyperbolic and parabolic trajectories in optical tweezers. 
The transition from mere optical manipulation towards optical control of matter, and further uses in other fields, \textit{e.g.} in quantum information protocols \cite{Leach2010p662,Jha2010p010501}, implies measuring the linear and orbital angular momenta, or compositions of these, of arbitrary light fields. Here, a simple schema to calculate momenta for scalar wave fields through the Fourier transform, based on the Hilbert space model for separation of variables of Helmholtz equation \cite{Boyer1976p35,Miller1984}, is presented for the benefit of advanced undergraduate and graduate students interested in the OAM of light and related topics. 

In the following section, a brief overview of the Helmholtz equation, its symmetries and separable solutions, is given; the separable solutions for Cartesian, circular-cylindrical and elliptic-cylindrical coordinate systems are shown as an example. In the third section, wave decompositions and integral transforms related to the Helmholtz equation are treated; Plancherel theorem, the keystone of the proposed schema, is presented in this section.
Fourier and Circular-cylindrical wave transforms of the separable solutions presented in section 2 are calculated for the sake of clarity. Next, the actual schema to calculate momenta through Fourier transforms is discussed and the mean value of the orbital angular momentum for Plane and Elliptical waves is calculated. 
Finally, a summary and plausible applications are discussed.

\section{Symmetries of the Hemlholtz equation} 
The Helmholtz differential equation, \begin{eqnarray} \left(\nabla^2 + k^2 \right) \Psi = 0, \end{eqnarray} where the symbol $\nabla^2$ is the Laplacian operator and $k$ is the magnitude of the wave vector $\mathbf{k} = k (\sin\vartheta \cos\varphi,\sin\vartheta \sin\varphi,\cos\vartheta)$, and $\vartheta$ and $\varphi$ are the Euler angles, is separable in eleven coordinate systems \cite{Eisenhart1934p427,Morse1953v1,Arscott1981p33}. 
Each of these coordinate system defines a three element set of operators that span a Lie algebra \cite{Boyer1976p35,Miller1984}. These operators fulfil the eigenvalue equation 
\begin{eqnarray} \label{eq:ConservedOperators} 
\mathcal{O}^{(i)}_{j} \Psi = o^{(i)}_{j} \Psi, 
\end{eqnarray} 
with $i=1, \ldots,11$, and $j = 1,2,3$. By defining a characteristic label set given by the three eigenvalues corresponding to a given coordinate system, $\Lambda_{i} = \left\{ o^{(i)}_{j} \right\}$, it is possible to label the set of functions that solve the reduced wave equation in the $i$-th coordinate system as $\{ \Psi_{\Lambda_{i}} \}$. 
A common element to all characteristic label sets is the square linear momentum that fulfils $\mathbf{p}^2 \Psi_{\Lambda_{i}} \equiv - \nabla^{2} \Psi_{\Lambda_{i}} = \mathbf{k}^2 \Psi_{\Lambda_{i}}$; \textit{i.e.}, the first element of each and every characteristic label set $\Lambda_{j}$ is the magnitude of the wave vector, $o^{(i)}_{1} = k$. 
The four coordinate systems with axial symmetry, Cartesian $(i=1)$, circular- $(i=2)$, elliptic- $(i=3)$, parabolic- $(i=4)$ cylindrical, share the $z$-component of the linear momentum as the second dynamical variable, $p_{z} \Psi_{\Lambda_{i}} = (k \cos\theta) \Psi_{\Lambda_{i}}$; \textit{i.e.}, for these four coordinate systems the second element of each of the characteristic label set is the Euler angle $\theta$, $o^{(i)}_{2} = \theta$ with $i=1,\ldots,4$. 
There are five coordinate systems that share the orbital angular momentum, $l_{z} \Psi_{\Lambda_{i}} = n \Psi_{\Lambda_{i}}$, with $\mathbf{l} = \mathbf{r} \times \mathbf{p}$; for these, the second element of each of the characteristic label set is the discrete orbital angular momentum number $n$, also known as topological charge. 
Finally, there are another thirteen operators that can be easily found in the literature \cite{Boyer1976p35,Miller1984}. 

In optics, there is a particular interest in propagation invariant scalar wave fields. There exist four fundamental separable families with this characteristic; for the sake of brevity, we will present just the cases for plane, circular- and elliptic-cylindrical waves. 

\subsection{Plane waves} 
In the Cartesian coordinate system, $\mathbf{x}=($x$, $y$, $z$)$, the separable solution of the Helmholtz differential equation are plane waves, \begin{eqnarray}
\Psi_{\Lambda_{1}} &=& \left( \sin \vartheta \right)^{1/2} ~e^{\imath k (x \sin\vartheta \cos\varphi + y \sin\vartheta \sin\varphi + z \cos\vartheta)}, \mbox{ } -\pi\le\varphi<\pi, ~ 0\le \vartheta <\pi, \nonumber \\
\end{eqnarray} 
The conserved operators are the linear momenta $p_y$ and $p_z$, \begin{eqnarray} 
p_{y} \Psi_{\Lambda_{1}} &=& (k \sin\vartheta \sin\varphi) \Psi_{\Lambda_{1}},\\ 
p_{z} \Psi_{\Lambda_{1}} &=& (k \cos\vartheta) \Psi_{\Lambda_{1}}, \end{eqnarray} 
yielding a characteristic label set given by the magnitude of the wave vector
and the two Euler angles, $\Lambda_{1} = \left\{ k, \vartheta, \varphi \right\}$. 

\subsection{Circular-cylindrical waves} 
For circular-cylindrical coordinates, $\mathbf{x} = (r\cos\phi, r\sin\phi, z)$, the separable solution of the Helmholtz differential equation are circular-cylinder waves, 
\begin{eqnarray}
\Psi_{\Lambda_{2}} &=& \imath^n \left( 2\pi\sin\vartheta\right)^{1/2} J_{n}(k\sin\vartheta ~r)e^{\imath (n\phi + k z \cos\vartheta)}, \mbox{ } 0\le \vartheta <\pi, ~n=0,\pm1,\pm2,\ldots \nonumber \\ 
\end{eqnarray} 
where the function $J_{n}(\cdot)$ is the Bessel function of the first kind.
These scalar wave fields, so called Bessel waves, conserve the linear momenta $p_{z}$ and the orbital angular momenta $l_{z}$, 
\begin{eqnarray} p_{z} \Psi_{\Lambda_{2}} &=& (k \cos\vartheta) \Psi_{\Lambda_{2}}, \\ 
l_z \Psi_{\Lambda_{2}} &=& n \Psi_{\Lambda_{2}}, 
\end{eqnarray} 
thus the characteristic label set is given by the magnitude of the wave vector, one of the Euler angles and the orbital angular momentum number, or topological charge, $\Lambda_{2} = \left\{ k, \vartheta, n \right\}$, in that order.

\subsection{Elliptic-cylindrical waves}
In elliptic-cylindrical coordinates ($x=f\cosh\xi\cos\eta$,$y=f\sinh\xi\sin\eta$, $z$), the separable solution of the Helmholtz wave equation are the so called even and odd Mathieu waves, in that order,
\begin{eqnarray} 
\Psi c_{\Lambda_{3}} &=& \left(\sin\vartheta\right)^{1/2} c_n \textrm{Ce}(\xi,
q) \textrm{ce}(\eta,q) e^{\imath k z \cos\vartheta}, \quad ~ 0\le \vartheta <\pi,~n=0,2,\ldots\\ 
\Psi s_{\Lambda_{3}} &=& \left(\sin\vartheta\right)^{1/2} s_n \textrm{Se}(\xi, q) \textrm{se}(\eta,q) e^{\imath k z \cos\vartheta}, \quad ~ 0\le \vartheta <\pi,~n=1,2,\ldots 
\end{eqnarray} 
where the normalization coefficients are given by the expressions
\begin{eqnarray} 
c_{2n}&=& \frac{1}{A_{0}^{(2n)}} \textrm{ce}_{2n}(0;q) \textrm{ce}_{2n}(\pi/2;q), \\
c_{2n+1}&=& - \frac{1}{\sqrt{q} A_{1}^{(2n+1)}} \textrm{ce}_{2n+1}(0;q) \left. \partial_{v}\textrm{ce}_{2n+1}(v;q)\right|_{v=\pi/2}, \\
s_{2n+1}&=& \frac{1}{\sqrt{q} B_{1}^{(2n+1)}} \left. \partial_{v}\textrm{se}_{2n+1}(v;q)\right|_{v=0} \textrm{se}_{2n+1}(\pi/2;q), \\ 
s_{2n+2}&=& \frac{1}{q B_{2}^{(2n+2)}} \left. \partial_{v}\textrm{se}_{2n+2}(v;q)\right|_{v=0} \left. \partial_{v}\textrm{se}_{2n+2}(v;q)\right|_{v=\pi/2},
\end{eqnarray} 
where $ce_n$ ($se_n$) and $Ce_n$ ($Se_n$) are the ordinary and extraordinary $n$-th order Mathieu functions of even (odd) parity, sometimes named angular and radial Mathieu functions, in that order. 
The separation constant is given by $ q = (f k \sin \vartheta / 2)^2$ with $2f$ being the inter-foci distance that defines a particular elliptic cylindrical coordinate system. 
Mathieu waves conserve, apart from the $z$-component of the linear momenta, a linear superposition of orbital angular momenta $l_{z}$ and linear angular momenta $p_{x}$ 
\begin{eqnarray} 
(l_z^2+f^2 p_x^2)\Psi c_{\Lambda_{3}}&=& a_{n} \Psi c^{(3)}_{\{k,\gamma,n\}},\\
(l_{z}^2+f^2 p_x^2)\Psi s_{\Lambda_{3}}&=& b_{n} \Psi s^{(3)}_{\{k,\gamma,n\}},
\end{eqnarray} 
such that the characteristic label set is given by the magnitude of the wave vector, one Euler angles related to the $z-$component of the wave vector and the discrete order number, $\Lambda_{3} = \left\{ k, \vartheta, n \right\}$ in that order. 
The even (odd) eigenvalue $a_{n}$ ($b_{n}$) is fixed for a given combination of $n$ and $q$. 

\section{Wave decompositions and integral transforms} 
Each one of the eleven sets of separable functions, $\{ \Psi_{\Lambda_{i}} \}$, that solve the reduced wave equation forms an orthonormal basis, \begin{eqnarray} \label{eq:Normalization} 
\int d\mathbf{x} ~ \bar{\Psi}_{\Lambda_{i}} (\mathbf{x}) \Psi_{\Lambda_{i}^{\prime}}(\mathbf{x}) \equiv \langle \Psi_{\Lambda_{i}} \vert \Psi_{\Lambda_{i}^{\prime}} \rangle = \tilde{\delta}(\Lambda_{i} - \Lambda_{i}^{\prime}), 
\end{eqnarray}
where the bracket, $\langle \cdot, \cdot \rangle$, notation is used as a shorthand for integration over configuration space, the upper bar notation, $\bar{\cdot}$, stands for complex conjugation and the notation $\tilde{\delta}(\cdot)$ stands for a combination of Dirac and Kronecker deltas particular to each coordinate system. 
A $\Psi_{\Lambda_{i}}$-wave decomposition of an arbitrary well behaved function $\Phi(\mathbf{x})$ can be defined for each of the eleven orthonormal basis,
\begin{eqnarray}
\label{eq:WaveDecomposition} \Phi(\mathbf{x}) &=& \int d\Lambda_{i} ~ \phi(\Lambda_{i}) \Psi_{\Lambda_{i}}(\mathbf{x}), 
\end{eqnarray} 
alongside an integral transform, 
\begin{eqnarray} \label{eq:IntegralTransform} 
\phi(\Lambda_{i}) &=& \int d\mathbf{x} ~ \bar{\Psi}_{\Lambda_{i}}(\mathbf{x}) \Phi(\mathbf{x}) \equiv \langle \Psi_{\Lambda_{i}} \vert \Phi \rangle.
\end{eqnarray} 
It is possible to show that each and every one of these eleven wave decompositions can be related to each other through their amplitudes, by substituting the corresponding wave decomposition, Eq.(\ref{eq:WaveDecomposition}), of $\bar{\Psi}_{\Lambda_{i}}(\mathbf{x})$ and $ \Phi(\mathbf{x})$ into the integral transform given by Eq.(\ref{eq:IntegralTransform}),
\begin{eqnarray} 
\phi(\Lambda_{i}) &=& \int d\Lambda_{j} ~\bar{\psi}_{\Lambda_{i}}(\Lambda_{j}) \phi(\Lambda_{j}), 
\end{eqnarray} 
In shorthand, 
\begin{eqnarray} 
\langle \Psi_{\Lambda_{i}} \vert \Phi \rangle &=& \langle \langle
\psi_{\Lambda_{i}} \vert \phi \rangle \rangle_{\Lambda_{j}},\label{eq:GenPlancherel} 
\end{eqnarray} 
where double angle-brackets, $\langle\langle \cdot \vert \cdot \rangle\rangle_{\Lambda_{j}}$, have been used to denote integration over the phase space defined by the characteristic label set $\Lambda_{j}$. 
The above expression, Eq.(\ref{eq:GenPlancherel}), is the well known Plancherel theorem \cite{Plancherel1910p289,Titchmarsh1925p279,Watson1933p156}; it is particularly useful when both the integral transforms related to a  $\Psi_{\Lambda_{j}}$-wave decomposition of $\Phi$ and $\Psi_{\Lambda_{i}}$ are known. 

Unitarity of the eleven Helmholtz-related integral transforms follows straightforward from Eq.(\ref{eq:GenPlancherel}) yielding Parceval theorem, 
\begin{eqnarray} 
\langle \Phi \vert \Phi \rangle = \langle \langle \phi \vert \phi \rangle \rangle_{\Lambda_{i}}. \label{eq:ITintermsFT} 
\end{eqnarray} 

\subsection{Plane wave decomposition and Fourier transform} 
Take as an example the plane wave decomposition 
\begin{eqnarray} 
\Phi (\mathbf{x}) = \int d{\Lambda_{1}} ~ \phi(\Lambda_{1}) (\sin \vartheta)^{1/2} e^{\imath \mathbf{k} \cdot \mathbf{x} } , 
\end{eqnarray} 
related to the Fourier transform 
\begin{eqnarray} 
\phi(\mathbf{k}) = \int d\mathbf{x}~ (\sin \vartheta)^{1/2} e^{-\imath \mathbf{k} \cdot \mathbf{x} } \Phi(\mathbf{x}). 
\end{eqnarray} 
The Fourier transform of a plane wave is given by the expression \begin{eqnarray}
\psi_{\Lambda_{1}^\prime}(\Lambda_{1})&=& \left( \sin\gamma \right)^{-1/2} \delta(\varphi-\varphi^{\prime})\delta(\vartheta-\vartheta^{\prime}),
\end{eqnarray} 
where the function $\delta(\cdot)$ is the Dirac delta function. A circular-cylindrical wave has a Fourier transform 
\begin{eqnarray}
\psi_{\Lambda_{2}^{\prime}}(\Lambda_{1})&=& \left( 2\pi\sin\vartheta^{\prime} \right)^{-1/2} e^{\imath n \varphi} \delta(\theta - \vartheta^{\prime}).
\end{eqnarray} 
For elliptic-cylindrical waves, Fourier transforms are 
\begin{eqnarray} 
\psi c_{\Lambda_{3}^{\prime}}(\Lambda_{1}) &=& \left( \pi \sin\vartheta \right)^{-1/2} \textrm{ce}_n(\varphi,q) \delta(\vartheta-\vartheta^{\prime}), \\
\psi s_{\Lambda_{3}^{\prime}}(\Lambda_{1}) &=& \left( \pi \sin\vartheta \right)^{-1/2} \textrm{se}_n(\varphi,q) \delta(\vartheta-\vartheta^{\prime}) .
\end{eqnarray}

\subsection{Circular-cylindrical wave decomposition and related
integral transform} 
Now, take the circular cylindrical wave decomposition, also known as Bessel wave decomposition,
\begin{eqnarray} 
\Phi (\mathbf{x}) &=& \sum_{n} \int_{0}^{\pi} \sin\vartheta d\vartheta ~ \phi(\Lambda_{2}) \imath^n (2\pi\sin\vartheta)^{1/2} J_{n}(k\sin\vartheta ~\rho)e^{\imath (n\varphi + k z \cos\vartheta)}. 
\end{eqnarray}
This circular wave decomposition defines a Bessel transform related to the Helmholtz equation; let us call it a Helmholtz-Bessel transform (not to be confused with the Fourier-Bessel, so called Hankel, transform) for the sake of having a short name for it,
\begin{eqnarray} 
\phi(\Lambda_{2}) &=& \langle \Psi_{\Lambda_{2}} \vert \Phi \rangle = \int d\mathbf{x}~ \Phi(\mathbf{x}) (-\imath)^n (2\pi\sin\vartheta)^{1/2} J_{n}(k\sin\vartheta ~\rho)e^{-\imath (n\varphi + k z \cos\vartheta)}. \nonumber \\ 
\end{eqnarray} 
Sometimes, it can be bothersome to calculate the value of this expression. If the Fourier transform $\phi(\Lambda_{1})$ of the function $\Phi(\mathbf{x})$ is known, it may be simpler to calculate the coefficients of the Bessel wave decomposition, by using Plancherel theorem, Eq.(\ref{eq:GenPlancherel}), as
\begin{eqnarray}
\phi(\Lambda_{2}) &=& (2\pi)^{-1/2} (\sin\vartheta)^{1/2} \int_{-\pi}^{\pi} d\varphi ~\phi(k,\vartheta,\varphi) e^{-\imath n \varphi}.  \label{eq:BTintermsFT} 
\end{eqnarray} 
For example, the Helmholtz-Bessel transform of a plane wave, described by the label set $\Lambda_{1}^{\prime}$, is given by the expression,
\begin{eqnarray}
\psi_{\Lambda_{1}^{\prime}}(\Lambda_{2})&=& \left( 2\pi\sin\vartheta \right)^{-1/2} e^{-\imath n \varphi} \delta(\vartheta - \vartheta^{\prime}),
\end{eqnarray}
while that of a Bessel wave is given by, 
\begin{eqnarray}
\psi_{\Lambda_{2}^{\prime}}(\Lambda_{2})&=& \left( \sin\vartheta \right)^{-1/2} \delta(\vartheta- \vartheta^{\prime}) \delta_{n,n^\prime} , 
\end{eqnarray} 
where the symbol $\delta_{\cdot,\cdot}$ stands for Kronecker delta. 
And that of a Mathieu wave is 
\begin{eqnarray} \label{eq:MTinBT} 
\psi c_{\Lambda_{3}^{\prime}}(\Lambda_{2}) &=& 2^{-1/2} A^{(n^{\prime})}_{2m + (n^{\prime} \mathrm{mod} 2)} (q^{\prime})\delta(\vartheta- \vartheta^{\prime}) \delta_{n,2m + (n^{\prime} \mathrm{mod} 2)} , \\ 
\psi s_{\Lambda_{3}^{\prime}}(\Lambda_{2}) &=& 2^{-1/2} B^{(n^{\prime})}_{2m + (n^{\prime} \mathrm{mod} 2)} (q^{\prime}) \delta(\vartheta- \vartheta^{\prime}) \delta_{n,2m + (n^{\prime} \mathrm{mod} 2)}, 
\end{eqnarray} 
with $m = 0, 1, 2, \ldots$ and the coefficients $A^{(i)}_{j}$ ( $B^{(i)}_{j}$) are the $j$-th coefficients of the cosine (sine) expansion of the ordinary Mathieu functions of $i$-th order. 

\section{Quantifying momenta through the Fourier transform} 
It is possible to calculate the mean value of one of the aforementioned sixteen conserved operators $\mathcal{O}^{(i)}_{j}$, as well as linear superpositions of them, for any given wave $\Psi$ using its $\Psi_{\Lambda_{i}}$-wave decomposition, 
\begin{eqnarray} 
\langle \Phi \vert \mathcal{O}^{(i)}_{j} \vert \Phi \rangle &=&\sum_{\Lambda_{i}} o^{(i)}_{j} \vert\langle \Psi_{\Lambda_{i}} \vert \Phi \rangle\vert^2 . 
\end{eqnarray} 
Via Plancherel theorem, Eq.(\ref{eq:GenPlancherel}), it is possible to give this expression in terms of Fourier transforms, which may be already available and provide an easier path to calculate the overlaps, 
\begin{eqnarray} \label{eq:MeanMomentaFT} 
\langle \Phi \vert \mathcal{O}^{(i)}_{j} \vert \Phi \rangle &=&\sum_{\Lambda_{i}} o^{(i)}_{j} \vert\langle\langle \psi_{\Lambda_{i}} \vert \phi \rangle\rangle_{\Lambda_{1}}\vert^2 , 
\end{eqnarray} 
where the phase space defined by the label set $\Lambda_{1} = \{k, \vartheta, \varphi \}$ yields 
\begin{eqnarray} 
\langle\langle \psi \vert \phi \rangle\rangle_{\Lambda_{1}} = \int_{0}^{\infty} dk \int_{-\pi}^{\pi} d\varphi \int_{0}^{\pi} \sin\vartheta d\vartheta ~ \bar{\psi}_{\Lambda_{i}}(\Lambda_{1}) \phi(\Lambda_{1}). 
\end{eqnarray} 

\subsection{Orbital angular momentum of plane waves } 
It is simple to calculate the well-known null mean value of the orbital angular momentum $l_{z}$ for a plane wave, $\Lambda_{1} = \left\{k, \vartheta, \varphi \right\}$, by using Eq.(\ref{eq:MeanMomentaFT}) and Eq.(\ref{eq:BTintermsFT}),
\begin{eqnarray} 
\langle \Psi_{\Lambda_{1}} \vert l_{z} \vert \Psi_{\Lambda_{1}} \rangle &=& \sum_{\tilde{n} = -\infty}^\infty \tilde{n} ~ \vert\langle\langle \psi_{\tilde{\Lambda}_{2}} \vert \psi_{\Lambda_{1}} \rangle\rangle_{\Lambda_{1}^{\prime}}\vert^2,\\ 
&=&  \sum_{\tilde{n} = -\infty}^\infty \tilde{n} \vert (\sin \vartheta)^{1/2} (2 \pi \sin \tilde{\vartheta} )^{-1/2} e^{\imath n \vartheta} \delta(\tilde{\vartheta} - \vartheta) \vert^2 \\
&=& 0. 
\end{eqnarray} 

\subsection{Orbital angular momentum of elliptic-cylindrical waves } 
Again, by use of Eq.(\ref{eq:MeanMomentaFT}) and Eq.(\ref{eq:MTinBT}), it is possible to calculate the mean value of the angular momenta for any given elliptic-cylindrical wave as 
\begin{eqnarray} 
\langle \Psi c_{\Lambda_{3}} \vert l_{z} \vert \Psi c_{\Lambda_{3}} \rangle &=& \sum_{\tilde{n} = -\infty}^\infty \tilde{n} ~ \vert\langle\langle \psi_{\tilde{\Lambda}_{2}} \vert \psi c_{\Lambda_{3}} \rangle\rangle_{\Lambda_{1}^{\prime}}\vert^2,\\
&=& \sum_{m = 0}^\infty \left[ m + (n \mathrm{mod} 2) /2 \right] \vert \left( \sin \vartheta / \sin \tilde{\vartheta} \right)^{1/2} A^{(n)}_{2m + (n \mathrm{mod} 2) } \delta(\tilde{\vartheta} - \vartheta) \vert^2 , \nonumber \\
&=& \sum_{m = 0}^\infty \left[ m + (n \mathrm{mod} 2) /2 \right] \vert A^{(n)}_{2m + (n \mathrm{mod} 2) } \delta(0) \vert^2 ,\\ \langle \Psi s_{\Lambda_{3}} \vert l_{z} \vert \Psi s_{\Lambda_{3}} \rangle &=& \sum_{\tilde{n} = -\infty}^\infty \tilde{n} ~ \vert\langle\langle \psi_{\tilde{\Lambda}_{2}} \vert \psi s_{\Lambda_{3}} \rangle\rangle_{\Lambda_{1}^{\prime}}\vert^2,\\ 
&=& \sum_{m = 0}^\infty \left[ m + (n \mathrm{mod} 2) /2 \right] \vert \left( \sin \vartheta / \sin \tilde{\vartheta} \right)^{1/2} B^{(n)}_{2m + (n \mathrm{mod} 2) }\delta(\tilde{\vartheta} - \vartheta) \vert^2 . \nonumber \\ 
&=& \sum_{m = 0}^\infty \left[ m + (n \mathrm{mod} 2) /2 \right] \vert B^{(n)}_{2m + (n \mathrm{mod} 2) }\delta(0) \vert^2 . 
\end{eqnarray} 
Note that the divergences in all these expressions are related to the fact that ideal scalar wave solutions to the differential Helmholtz equation have a divergent normalization, Eq.(\ref{eq:Normalization}).

\section{Final remarks} 
A schema to calculate momenta related to the symmetries of the Helmholtz differential equation has been presented. The procedure relies on the orthonormal bases and the sixteen conserved operators derived from the eleven coordinate systems where the reduced wave equation accepts a separable solution. 
It has been shown that all the integral transforms, or wave overlaps, yielding the coefficients for the wave decompositions related to the aforementioned orthonormal basis can be calculated through Fourier transforms via Plancherel theorem. These overlaps are particularly useful to calculate the mean values of the conserved operators, and linear superposition of them, constructed from the three linear and three angular momenta. 
As an example, the mean value of the orbital angular momentum of Elliptic-cylindrical waves, so called Mathieu waves, has been presented. With further efforts the schema can be extended from these ideal scalar wave fields to real world beams; by using the orthonormal bases given by Gauss, Bessel-, Mathieu- \cite{GutierrezVega2007p215} and Weber-Gauss \cite{RodriguezLara2010p327} beams it may be possible to implement computational procedures to calculate momenta of captured experimental images. 
An extension of these procedures to non-paraxial optical vector fields is plausible. Moreover, integral transforms and decompositions are of use for those in the quantum optics fields as they define decompositions between boson modes; \textit{e.g.} expressing creation/annihilation operators of Bessel
modes in terms of creation/annihilation operators of Plane wave
modes \cite{BialynickiBirula2006p342}. 

Furthermore, to give just an example from outside the field of optics, in quantum mechanics, the treatment of molecular torsional levels for water and similar molecules is related to Ince polynomials that are a solution to the Helmholtz differential equation in paraboloidal coordinates \cite{Roncaratti2010p716}. 
In these cases, the presented schema can help in calculating suitable combinations of linear and orbital angular momenta of the different torsional eigenstates of the molecules. 

\acknowledgments
The author is grateful to Prof. Dr. Eugenio Ley Koo who posed the interesting question on how to calculate momenta of arbitrary wave functions through wave decompositions and gave valuable comments and suggestions through the writing stage of this manuscript.

\section*{References}

\end{document}